\documentclass{gtmon_a}
\pdfoutput=1
\input{diagxy}

\proceedingstitle{Exotic homology manifolds (Oberwolfach 2003)}
\conferencestart{29 June 2003}
\conferenceend{5 July 2003}
\conferencename{Workshop on Exotic Homology Manifolds}
\conferencelocation{Oberwolfach Mathematics Institute, Oberwolfach,
Germany}

\editor{Frank Quinn}
\givenname{Frank}
\surname{Quinn}

\editor{Andrew Ranicki}
\givenname{Andrew}
\surname{Ranicki}

\title{The Bryant--Ferry--Mio--Weinberger construction of generalized manifolds}

\author{Friedrich Hegenbarth}
\givenname{Friedrich}
\surname{Hegenbarth}
\address{Department of Mathematics\\University of Milano\\\newline
Via C Saldini
50\\Milano\\Italy 02130}
\email{friedrich.hegenbarth@mat.unimi.it}
\urladdr{http://www.mat.unimi.it/}

\author{Du\v san Repov\v s}
\givenname{Du\v san}
\surname{Repov\v s}
\address{Institute for Mathematics, Physics and Mechanics\\University of
Ljubljana\\\newline
Jadranska 19\\Ljubljana\\Slovenia 1001}
\email{dusan.repovs@fmf.uni-lj.si}
\urladdr{http://pef.pef.uni-lj.si/~dusanr/index.htm}

\volumenumber{9}
\issuenumber{}
\publicationyear{2006}
\papernumber{3}
\startpage{17}
\endpage{32}

\doi{}
\MR{}
\Zbl{}


\received{7 July 2005}
\revised{}
\accepted{7 July 2005}
\published{22 April 2006}
\publishedonline{22 April 2006}
\proposed{}
\seconded{}
\corresponding{}
\version{}

\arxivreference{math/0608654} 

\makeatletter
\def\cnewtheorem#1[#2]#3{\newtheorem{#1}{#3}[section]
\expandafter\let\csname c@#1\endcsname\c@thm}

\newtheorem{thm}{Theorem}[section]
\cnewtheorem{lem}[thm]{Lemma}
\cnewtheorem{cor}[thm]{Corollary}
\newtheorem*{thm'}{\fullref{thm-three-five}$'$}
\theoremstyle{definition}
\cnewtheorem{defn}[thm]{Definition}
\newtheorem*{rem}{Remark}

\makeatother  

\renewcommand{\to}{\rightarrow}
\let\barrsquare\square
\let\square\undefined
\makeop{Id}
\makeop{pr}

\begin{document}

\begin{htmlabstract}
Following Bryant, Ferry, Mio and Weinberger we construct
generalized manifolds as limits of controlled sequences {p<sub>i</sub>:
X<sub>i</sub>&rarr; X<sub>i-1</sub> : i = 1,2,&hellip;}
of controlled Poincar&eacute; spaces. The basic ingredient is the
&epsilon;-&delta;&ndash;surgery sequence recently proved by Pedersen,
Quinn and Ranicki. Since one has to apply it not only in cases when the
target is a manifold, but a controlled Poincar&eacute; complex, we explain
this issue very roughly. Specifically,
it is applied in the inductive step to construct the desired controlled
homotopy equivalence p<sub>i+1</sub>: X<sub>i+1</sub>&rarr;X<sub>i</sub>.
Our main theorem requires a sufficiently controlled Poincar&eacute;
structure on X<sub>i</sub> (over X<sub>i-1</sub>). Our construction shows that this can
be achieved. In fact, the Poincar&eacute; structure of X<sub>i</sub> depends upon a
homotopy equivalence used to glue two manifold pieces together (the rest
is surgery theory leaving unaltered the Poincar&eacute; structure). It follows
from the &epsilon;-&delta;&ndash;surgery sequence (more precisely from
the Wall realization part) that this homotopy equivalence is sufficiently
well controlled.  In the final section we give additional explanation
why the limit space of the X<sub>i</sub>'s has no resolution.
\end{htmlabstract}

\begin{abstract}
Following Bryant, Ferry, Mio and Weinberger we construct
generalized manifolds as limits of controlled sequences $\{p_i\colon
X_i\longrightarrow X_{i-1} : i = 1,2,\ldots\}$
of controlled Poincar\'e spaces. The basic ingredient is the
$\varepsilon{-}\delta$--surgery sequence recently proved by Pedersen,
Quinn and Ranicki. Since one has to apply it not only in cases when the
target is a manifold, but a controlled Poincar\'e complex, we explain
this issue very roughly. Specifically,
it is applied in the inductive step to construct the desired controlled
homotopy equivalence $p_{i+1}\colon X_{i+1}\to X_i$.
Our main theorem requires a sufficiently controlled Poincar\'e
structure on $X_i$ (over $X_{i-1}$). Our construction shows that this can
be achieved. In fact, the Poincar\'e structure of $X_i$ depends upon a
homotopy equivalence used to glue two manifold pieces together (the rest
is surgery theory leaving unaltered the Poincar\'e structure). It follows
from the $\varepsilon{-}\delta$--surgery sequence (more precisely from
the Wall realization part) that this homotopy equivalence is sufficiently
well controlled.  In the final section we give additional explanation
why the limit space of the $X_i$'s has no resolution.
\end{abstract}

\begin{htmlabstract}
Following Bryant, Ferry, Mio and Weinberger we construct
generalized manifolds as limits of controlled sequences {p<sub>i</sub>:
X<sub>i</sub>&rarr; X<sub>i-1</sub> : i = 1,2,&hellip;}
of controlled Poincar&eacute; spaces. The basic ingredient is the
&epsilon;{-}&delta;&ndash;surgery sequence recently proved by Pedersen,
Quinn and Ranicki. Since one has to apply it not only in cases when the
target is a manifold, but a controlled Poincar&eacute; complex, we explain
this issue very roughly. Specifically,
it is applied in the inductive step to construct the desired controlled
homotopy equivalence p<sub>i+1</sub>:X<sub>i+1</sub>&rarr; X<sub>i</sub>.
Our main theorem requires a sufficiently controlled Poincar&eacute;
structure on X<sub>i</sub> (over X<sub>i-1</sub>). Our construction shows that this can
be achieved. In fact, the Poincar&eacute; structure of X<sub>i</sub> depends upon a
homotopy equivalence used to glue two manifold pieces together (the rest
is surgery theory leaving unaltered the Poincar&eacute; structure). It follows
from the &epsilon;{-}&delta;&ndash;surgery sequence (more precisely from
the Wall realization part) that this homotopy equivalence is sufficiently
well controlled.  In the final section we give additional explanation
why the limit space of the X<sub>i</sub>'s has no resolution.
\end{htmlabstract}

\begin{abstract}
Following Bryant, Ferry, Mio and Weinberger we construct
generalized manifolds as limits of controlled sequences $\bigl\{
X_i \smash{\stackrel{p_i}{\longrightarrow}} X_{i-1} : i = 1,2,\ldots\bigl\}$
of controlled Poincar\'e spaces. The basic ingredient is the
$\varepsilon{-}\delta$--surgery sequence recently proved by Pedersen,
Quinn and Ranicki. Since one has to apply it not only in cases when the
target is a manifold, but a controlled Poincar\'e complex, we explain
this issue very roughly. Specifically,
it is applied in the inductive step to construct the desired controlled
homotopy equivalence $p_{i+1} \co  X_{i+1} \to X_i$.
Our main theorem requires a sufficiently controlled Poincar\'e
structure on $X_i$ (over $X_{i-1}$). Our construction shows that this can
be achieved. In fact, the Poincar\'e structure of $X_i$ depends upon a
homotopy equivalence used to glue two manifold pieces together (the rest
is surgery theory leaving unaltered the Poincar\'e structure). It follows
from the $\varepsilon{-}\delta$--surgery sequence (more precisely from
the Wall realization part) that this homotopy equivalence is sufficiently
well controlled.  In the final section we give additional explanation
why the limit space of the $X_i$'s has no resolution.
\end{abstract}

\maketitle

\section{Preliminaries}
\label{sec-one}

A \emph{generalized $n$--dimensional manifold} $X$ is
characterized by the following two properties:
\begin{enumerate}
\item[(i)] $X$ is a Euclidean neighborhood retract (ENR); and
\item[(ii)] $X$ has the local homology (with integer coefficients)
of the Euclidean $n$--space $\mathbb{R}^n$, ie
$$H_* (X, X \setminus \{x\}) \cong H_* (\mathbb{R}^n, \mathbb{R}^n \setminus \{0\}).$$
\end{enumerate}
Since we deal here with locally compact separable
metric spaces of finite (covering) dimension, ENRs are the
same as ANRs.

Generalized manifolds are Poincar\'e spaces,
in particular they have the Spivak normal fibrations
$\nu_X$. The total space of $\nu_X$ is the boundary
of a regular neighborhood $N(X) \subset \mathbb{R}^L$ of an
embedding $X \subset \mathbb{R}^L$, for some large $L$. One can assume
that $N(X)$ is a mapping cylinder neighborhood
(see Lacher \cite[Corollary 11.2]{La}).

The global Poincar\'e duality of Poincar\'e spaces
does not imply the local homology condition (ii) above.
The local homology condition can be understood
as the ``controlled'' global Poincar\'e duality
(see Quinn \cite[p270]{Qu83}, and Bryant--Ferry--Mio--Weinberger\break
\hbox{\cite[Proposition~4.5]{BFMW}}).
More precisely, one has the following:

\begin{thm}
\label{thm-one-one}
Let $X$ be a compact
ANR Poincar\'e duality space of finite
(covering) dimension. Then $X$ is a
generalized manifold if and only if for every $\delta > 0$, $X$ is a $\delta$--Poincar\'e space (over $X$).
\end{thm}

The definition of the $\delta$--Poincar\'e property is given below.
The following basic fact about homology
manifolds was proved by Ferry and Pedersen \cite[Theorem 16.6]{FP}.

\begin{thm} 
\label{thm-one-two}
Let $X$ be an ANR homology manifold. Then $\nu_X$
has a canonical {\rm TOP} reduction.
\end{thm}

This statement is equivalent to existence of degree-one normal maps
$f\co  M^n \to X$, where $M^n$ is a (closed) topological $n$--manifold,
hence the structure set ${\cal S}^\text{TOP} (X)$ can be identified with
$[X, G/\text{TOP}]$.

Let us denote the $4$--periodic simply connected surgery spectrum
by $\mathbb{L}$ and let $\widehat {\mathbb{L}}$ be the connected covering of
$\mathbb{L}$.  There is a (canonical) map of spectra $\widehat {\mathbb{L}}
\to \mathbb{L}$ given by the action of $\widehat {\mathbb{L}}$ on $\mathbb{L}$.
Note that $\widehat {\mathbb{L}}_0$ is $G/\text{TOP}$.

If $M^n$ is a topological manifold there exists a
fundamental class $[M]_{\mathbb{L}} \in H_n (M; \mathbb{L}^{\bullet})$, where
$\mathbb{L}^{\bullet}$ is the symmetric surgery spectrum
(see Ranicki \cite[Chapters 13 and 16]{Ra}).

\begin{thm}
\label{thm-one-three}
If $M^n$ is a closed oriented topological $n$--manifold,
then the cap product with $[M]_{\mathbb{L}}$ defines a Poincar\'e
duality of $\mathbb{L}$--(co)homology
$$
H^p (M; \mathbb{L}) \overset\cong\longrightarrow H_{n-p} (M; \mathbb{L})
$$
and $\widehat {\mathbb{L}}$--(co)homology
$$
H^p (M; \widehat {\mathbb{L}}) \overset\cong\longrightarrow H_{n-p} (M; \widehat {\mathbb{L}}).
$$
\end{thm}

Since $H^0 (M; \mathbb{L}) = [M, \mathbb{Z} \times G/\text{TOP}]$ and
$H^0 (M; \widehat {\mathbb{L}}) = [M, G/\text{TOP}]$, we have
$$H_n (M; \mathbb{L}) = \mathbb{Z} \times H_n (M; \widehat {\mathbb{L}})$$
and
the map $\widehat {\mathbb{L}} \to \mathbb{L}$ has the property that
the image of
$$H_n (M; \widehat {\mathbb{L}}) \to H_n (M; \mathbb{L}) = \mathbb{Z} \times H_n (M;
\widehat {\mathbb{L}})$$
is $\{1\} \times H_n (M; \widehat {\mathbb{L}})$
(see Ranicki \cite[Appendix C]{Ra}).
Moreover, the action of $H^0 (M; \widehat {\mathbb{L}})$ on
$H^0 (M; \mathbb{L}) = \mathbb{Z} \times H^0 (M; \widehat {\mathbb{L}})$, induced
by the action of $\widehat {\mathbb{L}}$ on $\mathbb{L}$, preserves
the $\mathbb{Z}$--sectors.

If $X$ is a generalized $n$--manifold
we get similar results by using the
fundamental class $f_*([M]_{\mathbb{L}}) = [X]_{\mathbb{L}} \in H_n (X; \mathbb{L}^{\bullet})$,
where $f\co M \to X$ is the
canonical degree-one normal map. So
the composition map
$$
\Theta \co  [X, G/\text{TOP}] \to H_n (X; \widehat {\mathbb{L}}) \to H_n (X; \mathbb{L}) =
\mathbb{Z} \times H_n (X; \widehat {\mathbb{L}})
$$
has the property that $\operatorname{Im} \Theta$ belongs to
a single $\mathbb{Z}$--sector, denoted by $I(X) \in \mathbb{Z}$.

The following is the fundamental result of Quinn on resolutions
of generalized manifolds \cite{Qu87}.

\begin{thm} 
\label{thm-one-four}
Let $X$ be a generalized $n$--manifold, $n \ge 5$.
Then $X$ has a resolution if and only if $I(X) = 1$.
\end{thm}

\begin{rem}
The integer $I(X)$ is called the \emph{Quinn index}
of the generalized manifold $X$. Since the
action of $\widehat {\mathbb{L}}$ on $\mathbb{L}$ preserves the $\mathbb{Z}$--sectors,
arbitrary degree-one normal maps $g \co  N \to X$
can be used to calculate $I(X)$.
Alternatively, we can define $I(X)$ using the
fibration $\widehat {\mathbb{L}} \to \mathbb{L} \to \mathbb{K} (\mathbb{Z}, 0)$,
where $\mathbb{K} (\mathbb{Z}, ^.)$ is the Eilenberg--MacLane spectrum,
and define $I(X)$ as the image of (see Ranicki \cite[Chapter 25]{Ra}):
$$\{f \co  M \to X\} \in H_n (X; \mathbb{L}) \to H_n (X; \mathbb{K} (\mathbb{Z}, 0))
  = H_n (X; \mathbb{Z}) = \mathbb{Z}.$$
We assume that $X$ is oriented. Therefore $I(X)$ is also defined
for Poincar\'e complexes, as long as we have a degree-one normal
map $f \co  M \to X$, determining an element in $H_n (X; \mathbb{L})$. In
this case $I(X)$ is not a local index. In fact, for generalized
manifolds one has local $\mathbb{L}$--Poincar\'e duality using locally
finite chains, hence we can define $I({\cal U})$ for any open set
${\cal U} \subset X$. It is also easy to see that $I({\cal U}) =
I(X)$. On the algebraic side $I(X)$ is an invariant of the
controlled Poincar\'e duality type (see Ranicki \cite[p283]{Ra}).
\end{rem}


\section{Constructing generalized manifolds from controlled
se\-qu\-ences of Poincar\'e complexes}
\label{sec-two}

Beginning with a closed topological $n$--manifold $M^n$, for $n \ge 5$,
and $\sigma \in H_n (M; \mathbb{L})$, we shall construct a sequence of closed
Poincar\'e duality spaces $X_0, X_1, X_2, \ldots$, and maps $p_i \co  X_i
\to X_{i-1}$ and $p_0 \co  X_0 \to M$.

We assume that $M$ is a PL manifold,
or that $M$ has a cell structure. The $X_i$
are built by gluing manifolds along boundaries
with homotopy equivalences, and by doing some
surgeries outside the singular sets. Hence all the
$X_i$ have cell decompositions.

We can assume that the $X_i$ lie in a
(large enough) Euclidean space $\mathbb{R}^L$ which induces the metric
on $X_i$. So the cell chain complex $C_{\#} (X_i)$
can be considered as a geometric chain complex
over $X_{i-1}$ with respect to $p_i \co  X_i \to X_{i-1}$, ie the
distance between two cells of $X_i$ over $X_{i-1}$ is the
distance between the images of the centers
of these two cells in $X_{i-1}$. Let us denote the distance
function by $d$.

We now list five properties of the sequence
$\{(X_i, p_i)\}_i$, including some definitions
and comments.
For each $i \ge 0$ we choose positive real numbers
$\xi_i$ and $\eta_i$.

\begin{enumerate}
\item [(i)]
$p_i \co  X_i \to X_{i-1}$ and $p_0 \co  X_0 \to M$
are ${\cal U}V^1$--maps.
This means that for every $\varepsilon > 0$ and for all diagrams
%
%
$$\bfig
\barrsquare/->`^{ (}->`->`->/<600,400>[K_0`X_i`K`X_{i-1};
\alpha_0``p_i`\alpha]
\morphism/-->/<600,400>[\phantom{K}`\phantom{X_i};\wbar\alpha]
\efig$$
%
%
with $K$ a $2$--complex, $K_0 \subset K$ a subcomplex and
maps $\alpha_0,\alpha$, there is a map $\overline{\alpha}$ such that
$\overline{\alpha}\big\vert_{K_0} = \alpha_0$ and $d(p_i \circ \overline{\alpha}, \alpha) < \varepsilon$.
(This is also called ${\cal U}V^1 (\varepsilon)$ property.)

\item[(ii)] $X_i$ is an $\eta_i$--Poincar\'e complex over $X_{i-1}$, ie

  \begin{enumerate}
  \item[(a)] all cells of $X_{i-1}$ have diameter $< \eta_i$ over
  $X_{i-1}$; and

  \item[(b)] there is an $n$--cycle $c \in C_n(X_i)$ which induces an
  $\eta_i$--chain equivalence $\cap_c \co  C^{\#} (X_i) \to C_{n-\#} (X_i)$.
  \end{enumerate}

Equivalently, the diagonal $\Delta_{\#} (c) = \sum c'\otimes c'' \in
C_{\#}(X) \otimes C_{\#}(X)$ has the property that $d(c', c'') < \eta_i$
for all tensor products appearing in $\Delta_{\#} (c)$.

\item[(iii)] $p_i \co  X_i \to X_{i-1}$ is an $\xi_i$--homotopy
equivalence over $X_{i-2}$, for $i \ge 2$.
In other words, there exist an inverse $p'_i \co  X_{i-1} \to X_i$
and homotopies $h_i \co  p'_i \circ p_i \simeq \operatorname{Id}_{X_i}$
and $h'_i \co  p_i \circ p'_i \simeq \operatorname{Id}_{X_{i-1}}$
such that the tracks
$$\bigl\{(p_{i-1} \circ p_i \circ h_i) (x, t) : t \in
[0,1]\bigr\}\quad\text{and}\quad
\bigl\{(p_{i-1} \circ h'_i) (x', t) : t \in [0,1]\bigr\}$$
have diameter less than
$\xi_i$, for each $x \in X_i$ (respectively, $x' \in X_{i-1}$). Note that $p_0$
need not be a homotopy equivalence.

\item[(iv)] There is a regular neighborhood $W_0 \subset \mathbb{R}^L$
of $X_0$ such that $X_i \subset W_0$, for $i = 0, 1,\ldots$, and
retractions $r_i \co  W_0 \to X_i$, satisfying $d(r_i, r_{i-1}) < \xi_i$
in $\mathbb{R}^L$.

\item[(v)] There are ``thin'' regular neighborhoods $W_i \subset
\mathbb{R}^L$ with $\pi_i \co  W_i \to X_i$, where $W_i \subset \smash{\mathop
W \limits^\circ} _{i-1}$ such that $W_{i-1} \setminus \smash{\mathop W
\limits^\circ} _i$ is an $\xi_i$--thin $h$--cobordism with respect to
$r_i \co  W_0 \to X_i$.

Let $W = W_{i-1} \setminus \smash{\mathop W \limits^\circ} _i$. Then there
exist deformation retractions $r^0_t \co  W \to \partial_0W$ and $r^1_t
\co  W \to \partial_1W$ with tracks of size $< \xi_i$ over $X_{i-1}$,
ie the diameters of $\bigl\{(r_i \circ r^0_t) (w) : t \in [0,1]\bigr\}$ and
$\bigl\{(r_i \circ r^1_t) (w) : t \in [0,1]\bigr\}$ are smaller than $\xi_i$.
Moreover, we can choose $\eta_i$ and $\xi_i$ such that

\begin{enumerate}
\item[(a)] $\sum \eta_i < \infty$; and

\item[(b)] $W_{i-1} \setminus \smash{\mathop W \limits^\circ} _i$
has a $\delta_i$--product structure with $\sum \delta_i < \infty$, ie
there is a homeomorphism
$$W = W_{i-1} \setminus
{\mathop W \limits^\circ}_i \; {\stackrel{H}{\gets}} \;
\partial_0 W \times I$$
satisfying
$$\operatorname{diam} \bigl\{(r_i \circ H) (w,t) : t \in I\bigr\} < \delta_i,$$
for every $w \in \partial_0 W$.
\end{enumerate}
\end{enumerate}

The property (v)(b) above follows from the ``thin $h$--cobordism''
theorem (see the article \cite{Qu79} by Quinn). One can assume that
$\sum \xi_i < \infty$.
Let $X = \bigcap_i W_i$. We are going to show that
$X$ is a generalized manifold:

\begin{enumerate}
\item[(1)] The map $r = \varinjlim r_i \co  W_0 \to X$ is well--defined
and is a retraction, hence $X$ is an ANR.

\item[(2)] To show that $X$ is a generalized manifold we shall apply
the next two theorems.  They also imply \fullref{thm-one-one}
above.  The first one is due to Daverman and Husch \cite{DH}, but it is
already indicated in \cite{Qu79} (see the remark
after Theorem 3.3.2).
\end{enumerate}

\begin{thm} 
\label{thm-two-one}
Suppose that $M^n$ is a closed topological $n$--manifold,
$B$ is an ANR, and
$p \co  M \to B$ is proper and onto. Then
$B$ is a generalized manifold,
provided that $p$ is an approximate fibration.
\end{thm}

Approximate fibrations are characterized by the property that
for every $\varepsilon > 0$ and every diagram
%
%
$$\bfig
\barrsquare/->`^{ (}->`->`->/<600,400>[K\times\{0\}`M`K\times I`B;
H_0``p`h]
\morphism/-->/<600,400>[\phantom{K\times I}`\phantom{M};H]
\efig$$
%
%
where $K$ is a polyhedron,
there exists a lifting $H$ of $h$ such that
$d(p \circ H, h) < \varepsilon$. Here $d$ is a metric on $B$.
In other words, $p \co  M \to B$ has the $\varepsilon$--homotopy lifting
property for all $\varepsilon > 0$.

We apply \fullref{thm-two-one} to the map $\rho \co  \partial W_0 \to X$
defined as follows: Let $\rho \co  W_0 \to X$ be the map which
associates to $w \in W_0$ the endpoint $\rho (x) \in X$ following
the tracks defined by the ``thin'' product structures of the
$h$--cobordism when decomposing
$$W_0 = (W_0 \setminus {\mathop W \limits^\circ} _1) \cup (W_1 \setminus
{\mathop W \limits^\circ} _2) \cup \dots$$
The restriction to $\partial W_0$ will also be denoted by $\rho$.
By (v)(b) above, the map $\rho$ is well--defined and continuous.
We will show that it is an $\varepsilon$--approximate
fibration for all $\varepsilon > 0$.

The map $\rho \co  W_0 \to X$ is the limit of maps
$\rho_i \co  W_0 \to X_i$, where $\rho_i$ is the composition
given by the tracks
$(W_0 \setminus {\mathop W \limits^o} _1) \cup (W_1 \setminus {\mathop W \limits^o} _2) \cup
\dots \cup (W_{i-1} \setminus {\mathop W \limits^o} _i)$
followed by $\pi_i \co  W_i \to X_i$.
The second theorem is due to Bryant, Ferry, Mio and Weinberger
\cite[Proposition 4.5]{BFMW}.

\begin{thm} 
\label{thm-two-two}
Given $n$ and $B$, there exist $\varepsilon_0 > 0$ and $T > 0$
such that for every $0 < \varepsilon < \varepsilon_0$ the following holds:
If $X \mathop \to \limits^p B$ is an $\varepsilon$--Poincar\'e complex
with respect to the ${\cal U}V^1$--map $p$ and
$W \subset \mathbb{R}^L$ is a regular neighborhood of $X \subset \mathbb{R}^L$, ie
$\pi \co  W \to X$ is a neighborhood retraction,
then $\pi \big\vert_{\partial W} \co  \partial W \to X$ has the
$T\varepsilon$--lifting property, provided that
the codimension of $X$ in $\mathbb{R}^L$ is $\ge 3$.
\end{thm}

This is applied as follows:
Let $B \subset \mathbb{R}^L$ be a (small) regular neighborhood
of $X \subset \mathbb{R}^L$. Hence $X_k \subset W_k \subset B$ for sufficiently
large $k$.
It follows by property (ii) that $X_i$ is
an $\eta_i$--Poincar\'e complex over $X_i \mathop \rightarrow \limits^{p_i} X_{i-1} \subset B$,
hence (for $i$ sufficiently large) we get the following:

\begin{cor}
\label{cor-two-three}
$\rho_i \co  \partial W_0 \to X_i$ is a
$T\eta_i$--approximate fibration over $B$.
\end{cor}

\begin{proof}
By the theorem above, $\pi_i \co  \partial W_i \to X_i$
is a $T\eta_i$--approximate fibration over $B$,
hence so is $\rho_i \co  \partial W_0 \cong \partial W_i \to X_i$.
\end{proof}

It follows by construction that
$
\smash{\varprojlim_{p_i}} X_i = X \subset B,
$
and so we have, in the limit,
an approximate fibration
$\rho \co  \partial W_0 \to X$ over $\operatorname{Id} \co  X \to X$, ie
$X$ is a generalized manifold.
We will show in \fullref{sec-four} that $I(X)$
is determined by the $\mathbb{Z}$--sector of $\sigma \in H_n (M; \mathbb{L})$.


\section{Construction of the sequence of controlled Poincar\'e complexes}
\label{sec-three}

Before we begin with the construction we need
more fundamental results about controlled surgery
and approximations.

\subsection{$\varepsilon{-}\delta$ surgery theory}

We recall the main theorem of the article \cite{PQR} by Pedersen, Quinn
and Ranicki.
Let $B$ be a finite--dimensional compact
ANR, and $N^n$ a compact $n$--manifold
(possibly with nonempty boundary $\partial N$), where $n \ge 4$.
Then there exists an $\varepsilon_0 > 0$ such that
for every $0 < \varepsilon < \varepsilon_0$ there exist $\delta > 0$
with the following property:

If $p \co  N \to B$ is a ${\cal U}V^1 (\delta)$ map, then there exists
a controlled exact surgery sequence
\begin{equation}
\label{eq-three-one}
H_{n+1} (B; \mathbb{L}) \to {\cal S}_{\varepsilon, \delta} (N,p) \to [N,
\partial N; G/\text{TOP}, *]
\smash{\stackrel{\Theta}{\longrightarrow}} H_n (B; \mathbb{L}).
\end{equation}
The controlled structure set ${\cal S}_{\varepsilon, \delta} (N, p)$ is
defined as follows.
Elements of ${\cal S}_{\varepsilon, \delta} (N, p)$ are (equivalence) classes
of $(M, g)$, where $M$ is an $n$--manifold,
$g \co  M \to N$
is a $\delta$--homotopy equivalence over $B$ and
$g \big\vert_{\partial M} \co  \partial M \to \partial N$ is a homeomorphism.
The pair
$(M, g)$ is related to $(M', g')$ if there is a homeomorphism
$h \co  M \to M'$, such that the diagram
%
%
$$\bfig
\Vtriangle<400,400>[\partial M`\partial M'`\partial N;h`g`g']
\efig$$
%
%
commutes,
and $g' \circ h$ is $\varepsilon$--homotopic to $g$ over $B$.
Since $\varepsilon$ is fixed, this relation is
not transitive. It is part of the assertion
that it is actually an equivalence relation.
Then ${\cal S}_{\varepsilon, \delta} (N, p)$ is the set of equivalence
classes of pairs $(M, g)$.

As in the classical surgery theory, the map
\begin{equation}
\label{eq-three-two}
H_{n+1} (B; \mathbb{L}) \to {\cal S}_{\varepsilon, \delta} (N,p)
\end{equation}
is the controlled realization of surgery obstructions, and
\begin{equation}
\label{eq-three-three}
{\cal S}_{\varepsilon, \delta} (N,p) \to [N, \partial N; G/\text{TOP}, *]
\mathop \to \limits^\Theta H_n (B; \mathbb{L})
\end{equation}
is the actual (controlled) surgery part.
The following discussion will show that \eqref{eq-three-three} also holds for
controlled Poincar\'e spaces
(see \fullref{thm-three-five} below).
Moreover, $\delta$ is also of (arbitrary) small size,
provided that such is also $\varepsilon$.

To see this we will go through some of the
main points of the proof of \cite[Theorem~1]{PQR}.
For $\eta, \eta' > 0$ we denote by $L_n (B, \mathbb{Z}, \eta, \eta')$ the
set of highly $\eta$--connected $n$--dimensional quadratic
Poincar\'e complexes modulo highly $\eta'$--connected
algebraic cobordisms. Then there is a well--defined
obstruction map
$$
\Theta_{\eta} \co  [N, G/\text{TOP}] \to L_n (B, \mathbb{Z}, \eta, \eta')
$$
(for simplicity we shall assume that $\partial N = \emptyset$). If
$(f, b) \co  M^n \to N^n$ is a degree-one normal map one can
do controlled surgery to obtain a highly $\eta$--connected
normal map $(f', b') \co  {M'}^n \to N^n$ over $B$. If $N^n$ is a
manifold this can be done for every $\eta > 0$. If $N^n$ is a
Poincar\'e complex, it has to be $\eta$--controlled over $B$.
By \fullref{thm-one-one} above, this holds in particular for
generalized manifolds.

Given $\eta > 0$ there is an $\eta' > 0$ such that
if $(f', b')$ and $(f'', b'')$ are normally bordant, highly
$\eta$--connected, degree-one, normal maps, there
is then a highly $\eta'$--connected normal bordism
between them.
(Again this is true if $N$ is
an $\eta$--Poincar\'e complex over $B$.) This defines
$\Theta_{\eta}$.

To eventually complete surgeries in the
middle dimension we assume that the map $p \co  N \to B$
is ${\cal U}V^1$. Then one has the following
(see \cite[p243]{PQR}).
Given $\delta > 0$ there exists $\eta > 0$ such that if
$\Theta_{\eta} \bigl([f', b']\bigr) = 0$, then $(f', b')$ is normally cobordant
to a $\delta$--homotopy equivalence.
Moreover, if $(f'', b'')$ and $(f', b')$ are highly
$\eta$--connected degree-one normal maps being
normally cobordant, then there is a
highly connected $\eta'$--bordism between them
(ie for given $\eta$ there is such an $\eta'$).
Then controlled surgery produces a
controlled $h$--cobordism which gives
an $\varepsilon$--homotopy by the thin $h$--cobordism
theorem.
This defines an element of ${\cal S}_{\varepsilon, \delta} (N, p)$,
and shows the semi--exactness of the sequence
\begin{equation}
\label{eq-three-four}
{\cal S}_{\varepsilon, \delta} (N,p) \to [N, G/\text{TOP}]
\mathop \to \limits^{\Theta_{\eta}} L_n (B, \mathbb{Z}, \eta, \eta'),
\end{equation}
ie that ${\cal S}_{\varepsilon, \delta} (N,p)$ maps onto the
kernel of $\Theta_{\eta}$.
We note that semi--exactness also holds for
$\eta$--controlled Poincar\'e complexes over $B$.

One cannot expect the sequence \eqref{eq-three-four} to be exact, ie that
the composition map is zero, since passing
from topology to algebra one loses control.
As it was noted by Pedersen, Quinn and Ranicki \cite[p243]{PQR}, $\varepsilon$ and $\delta$ are determined
by the controlled Hurewicz and Whitehead theorems.
Exactness of \eqref{eq-three-four} will follow by the Squeezing
Lemma of Pedersen and Yamasaki \cite[Lemma~4]{PY}.

The proof of \eqref{eq-three-three} will be completed by showing
that the assembly map
$$
A \co  H_n (B; \mathbb{L}) \to L_n (B, \mathbb{Z}, \eta, \eta')
$$
is bijective for sufficiently small $\eta$.
This follows by splitting the controlled quadratic
Poincar\'e complexes (ie the elements of $L_n (B', \mathbb{Z}, \eta, \eta')$) into
small pieces over small simplices of $B$ (we assume
for simplicity that $B$ is triangulated). If $\delta$ is
given, and if we want a splitting where each piece
is $\delta$--controlled, we must start the subdivision with
a sufficiently small $\eta$--controlled quadratic Poincar\'e
complex (see the following Remark). This can be done
by \cite[Lemma~6]{PQR} (see also Yamasaki \cite[Lemma 2.5]{Y}).
Since $A \circ \Theta = \Theta_{\eta}$, we get \eqref{eq-three-three} from
\eqref{eq-three-four}.
The stability constant $\varepsilon_0$ is determined by
the largest $\eta$ for which $A$ is bijective.

\begin{rem}
Yamasaki has estimated the size of $\eta$ in the Splitting Lemma.
If one performs a splitting so that the two summands are
$\delta$--controlled, then one needs an $\eta$--controlled
algebraic quadratic Poincar\'e complex with $\eta$ of size $\delta
/{(an^k{+}b)}$, where $a,b,k$ depend on $X$ ($k$ is
conjectured to be 1), and $n$ is the length of the complex. Of
course, squeezing also follows from the bijectivity of $A$ for
small $\eta$, but the result \cite[Lemma~3]{PY} of Pedersen and Yamasaki
is somehow a clean
statement to apply (see \fullref{thm-three-five} below). We also
note that the bijectivity of $A$ is of course, independent of whether $N$
is a manifold or a Poincar\'e complex.
\end{rem}

\begin{thm}
\label{thm-three-five}
Suppose that $N \smash{\stackrel{p}{\longrightarrow}} B$ is a ${\cal U}V^1$ map.
Let $\delta > 0$ be given (sufficiently small, ie $\delta < \delta_0$
for some $\delta_0$).
Then there is $\eta > 0$ (small with respect to $\delta$),
such that if $N$ is an $\eta$--Poincar\'e complex over $B$,
and $(f, b) \co  M \to N$ is a degree-one normal map,
then $\Theta (f, b) = 0 \in H_n (B; \mathbb{L})$ if (and only if)
$(f, b)$ is normally bordant to a $\delta$--equivalence.
\end{thm}

The ``only if'' part is more delicate and
follows by \cite[Lemma~3]{PY}.
So let $f\co  M^n \to N^n$ be a $\delta$--equivalence
defining a quadratic $\eta_1$--Poincar\'e complex $C$ in
$L_n (B, \mathbb{Z}, \eta_1, \eta'_1)$ which is $\eta_1$--cobordant to zero
via $[N, G/\text{TOP}] \to L_n (B, \mathbb{Z}, \eta_1, \eta'_1)$.

Then $C$ is $\kappa \eta_1$--cobordant to an arbitrary small
quadratic Poincar\'e complex (ie to a quadratic
$\eta$--complex) which is $\kappa \eta'_1$--cobordant to zero, with
$\eta_1$ sufficiently small (ie $\eta$ sufficiently small).
In this case we can also assume that $A$
is bijective. This proves the ``only if'' part.

\fullref{thm-three-five} can also be stated as follows:

\begin{thm'}
\label{thm-three-five-dash}
Let $N$ be a sufficiently fine
$\eta$--Poincar\'e complex over a ${\cal U}V^1$--map
$p \co  N \to B$. Then there exist $\varepsilon > 0$ and $\delta > 0$,
both sufficiently small, such that the sequence
$$
{\cal S}_{\varepsilon, \delta} (N, p) \to [N, G/\text{TOP}] \to H_n (B; \mathbb{L})
$$
is exact. In particular, it holds for generalized manifolds.
\end{thm'}


\subsection{${\cal U}V^1$ approximation}

Here we recall the results \cite[Proposition~4.3, Theorem~4.4]{BFMW} of
Bryant, Ferry, Mio and Weinberger.

\begin{thm} 
\label{thm-three-six}
Suppose that $f \co  (M^n, \partial M) \to B$ is a continuous map
from a compact $n$--manifold with boundary
such that the homotopy fiber of $f$ is simply
connected. If $n \ge 5$ then $f$ is homotopic
to a ${\cal U}V^1$--map. In case that $f \big\vert_{\partial M}$ is
already ${\cal U}V^1$, the homotopy is relative $\partial M$.
\end{thm}

We state the second theorem in the form which
we will need.

\begin{thm}[Ferry {{\cite[Theorem~10.1]{Fe}}}]
\label{thm-three-seven}
Let $p \co  N^n \to B$ be a map from a compact
$n$--manifold into a polyhedron, where $n \ge 5$.  Then:
\begin{enumerate}
\item[\rm(i)]
Given $\varepsilon > 0$, there is a $\delta > 0$, such that
if $p$ is a ${\cal U}V^1 (\delta)$--map then $p$ is $\varepsilon$--homotopic
to a ${\cal U}V^1$--map.
\item[\rm(ii)]
Suppose that $p \co  N \to B$ is a ${\cal U}V^1$ map. Then for each
$\varepsilon > 0$ there is a $\delta > 0$ (depending on $p$ and $\varepsilon$)
such that if $f \co  M \to N$ is a $(\delta{-}1)$--connected
map (over $B$) from a compact manifold $M$
of dimension at least 5, then $f$ is $\varepsilon$--close over
$B$ to a ${\cal U}V^1$--map $g \co  M \to N$.
\end{enumerate}
\end{thm}


\subsection{Controlled gluing}

\begin{thm}[Bryant--Ferry--Mio--Weinberger {{\cite[Proposition~4.6]{BFMW}}}]
\label{thm-three-eight}
Let $(M_1, \partial M_1)$ and $(M_2, \partial M_2)$
be (orientable) manifolds and $p_i \co  M_i \to B$
be ${\cal U}V^1$--maps. Then there exist $\varepsilon_0 > 0$ and $T > 0$
such that, for $0 < \varepsilon \le \varepsilon_0$ and
$h \co  \partial M_1 \to \partial M_2$ an (orientation preserving)
$\varepsilon$--equivalence, $M_1 \mathop \cup_h M_2$ is
a $T\varepsilon$--Poincar\'e complex over $B$.
\end{thm}


\subsection{Approximation of retractions}

\begin{thm}[Bryant--Ferry--Mio--Weinberger {{\cite[Proposition~4.10]{BFMW}}}]
\label{thm-three-nine}
Let $X$ and $Y$ be finite polyhedra.
Suppose that $V$ is a regular neighborhood of $X$
with $\dim V \ge 2\dim Y + 1$ and $r \co  V \to X$ is
a retraction.
If $f \co  Y \to X$ is an $\varepsilon$--equivalence
with respect to $p \co  X \to B$, then there exists
an embedding $i \co  Y \to V$ and a retraction $s \co  V \to i(Y)$ with
$d (p \circ r, p \circ s) < 2\varepsilon$.
\end{thm}

We now begin with the construction.
Let $M^n$ be a closed oriented (topological) manifold
of dimension $n \ge 6$. Let
$\sigma \in H_n (M; \mathbb{L})$ be fixed. Moreover, we assume
that $M$ is equipped with a simplicial structure.
Then let $M = B \cup_D C$
be such that $B$ is a regular neighborhood of
the $2$--skeleton, $D = \partial B$ is its boundary and
$C$ is the closure of the complement of $B$.
So $D = \partial C = B \cap C$ is of dimension $\ge 5$.

By \fullref{thm-three-six} above we can replace $(B, D) \subset M$ and $(C, D)
\subset M$, by ${\cal U}V^1$--maps $j \co  (B, D) \to M$ and $j \co
(C, D) \to M$, and realize $\sigma$ according to $H_n (M; \mathbb{L})
\to {\cal S}_{\varepsilon, \delta} (D, j)$ by a degree-one normal map
$F_{\sigma} \co  V \to D \times I$ with $\partial_0 V = D$, $\partial_1
V = D'$, $F_{\sigma} \big\vert_{\partial_0 V} = \operatorname{Id}$
and $f_{\sigma = }F_{\sigma} \big\vert_{\partial_1 V} \co  D' \to D$
a $\delta$--equivalence over $M$.

We then define $X_0 = B \cup_{f_{\sigma}} - V \cup_{\operatorname{Id}} C$, where
$-V$ is the cobordism $V$ turned upside down.
We use the map
$-F_{\sigma} \cup \operatorname{Id} \co  {-}V \cup_{\operatorname{Id}} C \to
D \times I \cup C \cong C$
to extend $j$ to a map
$p_0 \co  X_0 \to M$.

The Wall realization $V \to D \times I$ is such
that $V$ is a cobordism built from $D$ by adding
high--dimensional handles (similarly
beginning with $D'$). Therefore $p_0$ is a ${\cal U}V^1$ map:
If $(K,L)$ is a simplicial pair with $K$ a
$2$--complex, and if there is given
a diagram
%
%
$$\bfig
\barrsquare/->`^{ (}->`->`->/<600,400>[L`X_0`K`M;\alpha_0``p_0`\alpha]
\efig$$
%
%
then we first move (by an arbitrary small approximation)
$\alpha$ and $\alpha_0$ into $B$ by general position arguments.
Then one uses the ${\cal U}V^1$--property of $j \co  B \to M$.
By \fullref{thm-three-eight}, $X_0$ is a $T\delta$--Poincar\'e complex over $M$.
Note that we can choose $\delta$ as small as we want,
hence we get an $\eta_0$--Poincar\'e complex for a
prescribed $\eta_0$.
This completes the first step.

To continue
we define a manifold $M^n_0$ and a degree-one normal map
$g_0 \co  M^n_0 \to X_0$ by
$$
M_0 = B \cup_{\operatorname{Id}} V
\cup_{\operatorname{Id}} -V
\cup_{\operatorname{Id}} C \to
B \cup_{\operatorname{Id}} D \times I
\cup_{f_{\sigma}} -V
\cup_{\operatorname{Id}} C \cong X_0,
$$
using $F_{\sigma} \cup \operatorname{Id}\co  
V \cup_{\operatorname{Id}} -V \to
D \times I \cup_{f_{\sigma}} -V$.
By construction it has a controlled surgery
obstruction $\sigma \in H_n (M; \mathbb{L})$.

Moreover, there is $\overline{\sigma} \in H_n (X_0;\mathbb{L})$
with $p_{0*} (\overline{\sigma}) = \sigma$. This can be seen from
the diagram
%
%
$$\bfig
\barrsquare/->`<-``<-/<800,400>[{H_n(M_0;\mathbb{L})}`{H_n(X_0;\mathbb{L})}`
  {H^0(M_0;\mathbb{L})}`{H^0(X_0;\mathbb{L})};
  {g_{0}}_*`\cong``g_0^*]
\barrsquare(800,0)/->``<-`<-/<800,400>[\phantom{H_n(X_0;\mathbb{L})}`
  {H_n(M;\mathbb{L})}`
  \phantom{H^0(X_0;\mathbb{L})}`
  {H^0(M;\mathbb{L})};
  {p_0}_*``\cong`p_0^*]
\efig$$
%
%
The vertical isomorphisms are Poincar\'e dualities.
Since $p_0$ is a ${\cal U}V^1$ map, $\overline{\sigma}$ belongs to the same
$\mathbb{Z}$--sector as $\sigma$.
We will again denote $\overline{\sigma}$ by $\sigma$.

We construct $p_1 \co  X_1 \to X_0$ as above:
Let $M_0 = B_1 \cup_{D_1} C_1$, let $B_1$ be a regular neighborhood
of the $2$--skeleton (as fine as we want), let
$C_1$ be the closure of the complement and let
$D_1 = C_1 \cap B_1 = \partial C_1 = \partial B_1$, and $g_0 \co  D_1 \to
X_0$ be a ${\cal U}V^1$ map.
Then we realize $\sigma \in H_n (X_0; \mathbb{L}) \to {\cal S}_{\varepsilon_1, \delta_1} (D_1, g_0)$
by $F_{1, \sigma} \co  V_1 \to D_1 \times I$ with
$\partial_0 V_1 = D_1$, $\partial_1 V_1 = D'_1$,
$F_{1,\sigma} \big\vert_{\partial_0 V_1} = \operatorname{Id}$ and
$f_{1,\sigma} = F_{1,\sigma} \big\vert_{\partial_1 V_1} \co  D'_1 \to D_1$
a $\delta_1$--equivalence
over $X_0$.

We define $p'_1 \co  X'_1 \to X_0$ by
$$X'_1 = B_1 \mathop \cup \limits_{f_{1,\sigma}} -V_1 \mathop \cup
\limits_{\operatorname{Id}} C_1 \mathop \rightarrow \limits^{f'_1}
M_0 \cong B_1 \mathop \cup \limits_{\operatorname{Id}}
D_1 \times I \mathop \cup \limits_{\operatorname{Id}} C_1,$$
using $-F_{1,\sigma} \co  {-}V_1 \to D_1 \times I$, and then
$p'_1 = g_0 \circ f'_1 \co  X'_1 \to M_0 \to X_0$.

We now observe that
\begin{enumerate}
\item[(i)] by \fullref{thm-three-eight}, $X'_1$ is a
$T_1\delta_1$--Poincar\'e complex over $X_0$; and
\item[(ii)] $p'_1$ is a degree-one normal map with controlled surgery
obstruction
$$-p_{0*} (\overline{\sigma}) + \sigma = 0 \in H_n (M; \mathbb{L}).$$
\end{enumerate}

Let $\xi_1 > 0$ be given. We now apply
\fullref{thm-three-five} to produce a $\xi_1$--homotopy
equivalence by surgeries outside the
singular set (note that the surgeries which have
to be done are in the manifold part of $X'_1$).
For this we need a sufficiently small $\eta_0$--Poincar\'e
structure on $X_0$. However, this can
be achieved as noted above.
This finishes the second step.

We now proceed by induction.
What we need for the third step in order to
produce $p_2 \co  X_2 \to X_1$ is
\begin{enumerate}
\item[(i)]
a degree-one normal map $g_1 \co  M_1 \to X_1$
with controlled surgery obstruction
$\sigma \in H_n (X_0; \mathbb{L})$; and
\item[(ii)] $\overline{\sigma} \in H_n (X_1; \mathbb{L})$ with $p_{1*} (\overline{\sigma}) = \sigma$,
in the same $\mathbb{Z}$--sector as $\sigma \in H_n (X_0; \mathbb{L})$.
\end{enumerate}

One can get $g_1 \co  M_1 \to X_1$ as follows:
Consider $g'_1 \co  M'_1 \to X'_1$, where
$$
M'_1 = B_1 \cup_{\operatorname{Id}} V_1 \cup_{\operatorname{Id}}
-V_1 \cup_{\operatorname{Id}} C_1 \to
B_1 \cup_{\operatorname{Id}}
D_1 \times I \cup_{f_{1,\sigma}}
-V_1 \cup_{\operatorname{Id}} C_1 \cong X'_1
$$
is induced by $F_{1,\sigma} \co  V_1 \to D_1 \times I$ and the identity.
The map $g'_1$ is a degree-one normal map. Then one performs
the same surgeries on $g'_1$ as one has performed on $p'_1 \co  X'_1 \to X_0$
to obtain $X_1$. This produces the desired $g_1$.
For (ii) we note that $p_{1*}$ is a
bijective map preserving the $\mathbb{Z}$--sectors
(since $p_1$ is ${\cal U}V^1$).

So we have obtained the sequence of controlled
Poincar\'e spaces $p_i \co  X_i \to X_{i-1}$ and
$p_0 \co  X_0 \to M$ with degree-one normal
maps $g_i \co  M_i \to X_i$ and controlled
surgery obstructions $\sigma \in H_n (X_{i-1}; \mathbb{L})$.
The properties (iv) and (v) of \fullref{sec-two}
now follow by the thin $h$--cobordism
theorem and approximation of retraction.


\section{Nonresolvability, the DDP property and existence of generalized manifolds}
\label{sec-four}

\subsection{Nonresolvability}

At the beginning of the construction we have $\sigma \in H_n (M; \mathbb{L})$, where
$M$ is a closed (oriented) $n$--manifold with $n \ge 6$.
For each $m$ we constructed degree-one normal maps
$g_m \co  M_m \to X_m$ over $p_m \co  X_m \to X_{m-1}$, with
controlled surgery obstructions $\sigma_m \in H_n (X_{m-1}; \mathbb{L})$,
$p_{0*} (\sigma_1) = \sigma$, $p_{m*} (\sigma_{m+1}) = \sigma_m$, and
all $\sigma_m$ belong to the same $\mathbb{Z}$--sector as $\sigma$.
So we will call all of them $\sigma$.

We consider the normal map $g_m \co  M_m \to X_m$
as a controlled normal map over the identity map
$\operatorname{Id} \co  X_m \to X_m$,
and over $q_m \co  X_m \subset W_m
\smash{\stackrel{\rho}{\longrightarrow}} X$
(see \fullref{sec-two}). Since $\rho \big\vert_{\partial W_m}$ is an approximate
fibration and $d(r_i, r_{i-1}) < \xi_i$ and
$\sum^{\infty}_{i = m + 1} \xi_i < \varepsilon$,
for large $m$, we can assume that
$q_m$ is ${\cal U}V^1(\delta)$ for large $m$, so
$(q_m)_* \co  H_n (X_m;\mathbb{L}) \to H_n (X; \mathbb{L})$ maps $\sigma$ to
$(q_m)_* (\sigma) = \sigma'$, being in the same $\mathbb{Z}$--sector as $\sigma$.
The map $(q_m)_*$ is a bijective, and we denote $\sigma'$ by $\sigma$.
In other words, we have a surgery problem
%
%
$$\bfig
\barrsquare/->``->`/<600,400>[M_m`X_m``X;g_m``q_m`]
\efig$$
%
%
over $X$, 
with controlled surgery obstruction $\sigma \in H_n (X; \mathbb{L})$.
Our goal is to consider the surgery problem
%
%
$$\bfig
\barrsquare/->``->`/<600,400>[M_m`X_m``X;q_m \circ g_m``\Id`]
\efig$$
%
%
over $\operatorname{Id} \co  X \to X$,
and prove that $\sigma \in H_n (X; \mathbb{L})$ is its
controlled surgery obstruction.

Observe that $q_m$ is a $\delta$--homotopy equivalence
over $\operatorname{Id} \co  X \to X$ if $m$ is sufficiently large
(for a given $\delta$).

Let ${\cal N} (X) \cong [X, G/\text{TOP}]$ be the normal cobordism classes
of degree-one normal maps of $X$, and
let $HE_{\delta} (X)$ be the set of $\delta$--homotopy
equivalences of $X$ over $\operatorname{Id} \co  X \to X$.
Our claim will follow from the following lemma.

\begin{lem}
\label{lem-four-one}
Let $HE_{\delta'} (X) \times {\cal N} (X) \mathop \to \limits^{\mu} {\cal N} (X)$
be the action map, ie $\mu (h, f) = h \circ f$. Then
for sufficiently small $\delta' > 0$, the diagram
%
%
$$\bfig
  \barrsquare<1000,400>[HE_{\delta'}(X)\times\mathcal{N}(X)`\mathcal{N}(X)`
    \mathcal{N}(X)`{H_n(X;\mathbb{L})};
    \mu`\pr`\Theta`\Theta]
  \efig$$
%
%
commutes.
\end{lem}

\begin{proof}
This follows from \fullref{thm-three-five-dash} since
$HE_{\delta'} (X) \times {\cal S}_{\varepsilon'', \delta''} (X, \operatorname{Id}) \to
{\cal S}_{\varepsilon, \delta} (X, \operatorname{Id})$
for sufficiently small $\delta'$ and $\delta''$.
\end{proof}

We apply this lemma to the map
$HE_{\delta} (X_m, X) \times {\cal N} (X_m) \to {\cal N} (X)$, which sends
$(h, g)$ to $h \circ g$, where $HE_{\delta} (X_m, X)$ are
the $\delta$--homotopy equivalences $X_m \to X$ over $\operatorname{Id}_X$.
Let $\psi_m \co  X \to X_m$ be a controlled inverse of $q_m$.
Then $\psi_m$ induces
$$\psi_{m*} \co  HE_{\varepsilon} (X_m, X) \to HE_{\delta} (X),$$
where $\delta$ is some multiple of $\varepsilon$.
One can then write the following
commutative diagram (for sufficiently
small $\delta$).
%
%
$$\bfig
  \morphism(0,400)|l|<0,400>[
  {HE_{\varepsilon}(X_m,X)\times\mathcal{N}(X_m)}`
  {HE_{\varepsilon}(X_m,X)\times H_n(X_m;\mathbb{L})};
  \Id\times\Theta]

  \morphism(0,400)|l|<0,-400>[
  \phantom{HE_{\varepsilon}(X_m,X)\times\mathcal{N}(X_m)}`
  HE_{\delta}(X)\times\mathcal{N}(X);
  (\psi_m)_* \times (q_m)_*]

  \morphism(0,400)<1000,0>[
  \phantom{HE_{\varepsilon}(X_m,X)\times\mathcal{N}(X_m)}`
  \mathcal{N}(X);]

  \morphism(1000,400)|b|<600,0>[
  \phantom{\mathcal{N}(X)}`{H_n(X;\mathbb{L})};\Theta]

  \morphism(0,800)<1600,-400>[
  \phantom{HE_{\varepsilon}(X_m,X)\times H_n(X_m;\mathbb{L})}`
  \phantom{H_n(X;\mathbb{L})};]

  \morphism(0,0)|b|<1000,400>[
  \phantom{HE_{\delta}(X)\times\mathcal{N}(X)}`
  \phantom{\mathcal{N}(X)};
  \mu]

  \morphism(0,0)|b|<1000,0>[
  \phantom{HE_{\delta}(X)\times\mathcal{N}(X)}`\mathcal{N}(X);\pr]

  \morphism(1000,0)|b|<600,400>[
  \phantom{\mathcal{N}(X)}`
  \phantom{H_n(X;\mathbb{L})};
  \Theta]

  \efig$$
%
%
with $HE_{\varepsilon} (X_m,X) \times \mathcal{N}(X_m) \to H_n (X; \mathbb{L})$
given by $(h, \tau) \to h_*(\tau)$.

It follows from this that for large enough $m$,
$q_m \circ g_m \co  M_m \to X$ has controlled surgery
obstruction $\sigma \in H_n (X; \mathbb{L})$.
Hence we get non--resolvable generalized
manifolds if the $\mathbb{Z}$--sector of $\sigma$ is $\ne 1$.


\subsection{The DDP Property}

The construction allows one to get the DDP property for $X$
(see \cite[Section~8]{BFMW}). Roughly speaking, this can be
seen as follows. The first step in the
construction is to glue a highly connected
cobordism $V$
into a manifold
$M$ of dimension $n \ge 6$, in between the regular
neighborhood of the $2$--skeleton.

The result is a space which has the DDP.
The other constructions are surgery on middle--dimensional
spheres, which also preserves
the DDP. But since we have to take
the limit of the $X_m$'s, one must do it
more carefully (see \cite[Definition~8.1]{BFMW}):

\begin{defn}
\label{defn-four-two}
Given $\varepsilon > 0$ and $\delta > 0$, we say that a space $Y$ has
the $(\varepsilon, \delta)$--DDP if for each pair of maps $f, g \co  D^2 \to Y$
there exist maps $\overline{f}, \overline{g} \co  D^2 \to Y$ such that
$d (\overline{f} (D^2), \overline{g} (D^2)) > \delta$, $d(f, \overline{f}) < \varepsilon$ and
$d (g, \overline{g}) < \varepsilon$.
\end{defn}

\begin{lem}
\label{lem-four-three}
$\{ X_m \}$ have the $(\varepsilon, \delta)$--DDP
for some $\varepsilon > \delta > 0$.
\end{lem}

\begin{proof}
The manifolds $M^n_m$, for $n \ge 6$, have
the $(\varepsilon, \delta)$--DDP for all $\varepsilon$ and $\delta$. In fact,
one can choose a sufficiently fine triangulation, such that
any $f \co  D^2 \to M$ can be placed by arbitrary
small moves into the $2$--skeleton or into
the dual $(n{-}3)$--skeleton. Then $\delta$ is
the distance between these skeleta.
The remarks above show that the $X_m$
have the $(\varepsilon, \delta)$--DDP for some $\varepsilon$ and $\delta$.
\end{proof}

It can then be shown that $X = \varprojlim X_i$
has the $(2\varepsilon, \delta / 2)$--DDP
(see \cite[Proposition~8.4]{BFMW}).


\subsection{Special cases}

\begin{enumerate}
\item[(i)]
Let $M^n$ and $\sigma \in H_n (M; \mathbb{L})$ be given as above.
The first case which can occur is that
$\sigma$ goes to zero
under the assembly map $A \co  H_n (M; \mathbb{L}) \to L_n (\pi_1 M)$.
Then we can do surgery
on the normal maps $F_{\sigma} \co  V \to D \times I$,
$F_{1,\sigma} \co  V_1 \to D_1 \times I$ and so on, to replace
them by products. In this case the
generalized manifold $X$ is homotopy equivalent
to $M$.

\item[(ii)]
Suppose that $A$ is injective (or is an isomorphism).
Then $X$ cannot be homotopy equivalent to
any manifold, if the $\mathbb{Z}$--sector of $\sigma$ is $\ne 1$.
Suppose that $N^n \to X$ were a homotopy equivalence. It
determines an element in $[X, G/\text{TOP}]$ which
must map to $(1,0) \in H_n (X; \mathbb{L})$, because
its surgery obstruction in $L_n (\pi_1X)$ is zero and
$A$ is injective. This contradicts our assumption
that the index of $X$ is not equal to 1.
Examples of this type are given by the $n$--torus $M^n = T^n$.
\end{enumerate}


\subsection{Acknowledgements}

The second author was supported in part by the
Ministry for Higher Education, Science and Technology of the
Republic of Slovenia research grant. We thank M Yamasaki for helpful
insights and the referee for the comments.


\bibliographystyle{gtart}
\bibliography{link}

\end{document}